\documentclass[nosumlimits,twoside]{amsart}
\usepackage[cp850]{inputenc}
\usepackage[bookmarksnumbered,plainpages,hyperindex]{hyperref}
 \input{xy}
 \xyoption{all}
\newtheorem{theorem}{\sc Theorem}[section]
\newtheorem{proposition}[theorem]{\sc Proposition}
\newtheorem{notation}[theorem]{\sc Notation}
\newtheorem{lemma}[theorem]{\sc Lemma}
\newtheorem{corollary}[theorem]{\sc Corollary}
\theoremstyle{definition}
\newtheorem{definition}[theorem]{\sc Definition}

\theoremstyle{remark}

\newtheorem{claim}[theorem]{}

\def\ot{\otimes}

\def\K{\textrm{Ker}}
\def\C{\textrm{Coker}}
\def\Id{\textrm{Id}}
\def\M{\mathcal{M}}

\def\w{\wedge_E }
\begin{document}
\pagestyle{headings}
\title{The Heyneman-Radford Theorem for Monoidal Categories}
\author{A. Ardizzoni}
\address{University of Ferrara, Department of Mathematics, Via Machiavelli
35, Ferrara, I-44100, Italy} \email{alessandro.ardizzoni@unife.it}
\urladdr{http://www.unife.it/utenti/alessandro.ardizzoni}
\thanks{This paper was written while A. Ardizzoni was member of G.N.S.A.G.A. with partial financial support
from M.I.U.R..}
\begin{abstract}
We prove Heyneman-Radford Theorem in the framework of Monoidal
Categories.
\end{abstract} \maketitle
\section*{Introduction} \markboth{\sc{A. ARDIZZONI}} {\sc{The Heyneman-Radford Theorem for Monoidal
Categories}}

This paper is devoted to the proof of the Heyneman-Radford Theorem
for Monoidal Categories. The original Heyneman-Radford's Theorem
(see \cite[Proposition 2.4.2]{HR} or \cite[Theorem 5.3.1, page
65]{Mo}) is a very useful tool in classical Hopf algebra theory.
We also point out that our proof is pretty different from the
classical one and hence might be of some interest even in the
classical case.\newline \medskip

\textbf{Notations.} \ Let $[(X,i_{X})]$ be a subobject of an
object $E$ in an abelian category $\M,$ where
$i_{X}=i_{X}^{E}:X\hookrightarrow E$ is a monomorphism and
$[(X,i_{X})]$ is the associated equivalence class. By abuse of
language, we will say that $(X,i_{X})$ is a subobject of $E$ and
we will write $(X,i_{X})=(Y,i_{Y})$ to mean that $(Y,i_{Y}) \in
[(X,i_{X})]$. The same convention applies to cokernels. If
$(X,i_{X})$ is a subobject of $E$ then we will write
$(E/X,p_X)=\C(i_X)$, where $p_{X}=p_{X}^{E}:E\rightarrow
E/X$.\newline
Let $(X_{1},i_{X_{1}}^{Y_{1}})$ be a subobject of $Y_{1}$ and let $%
(X_{2},i_{X_{2}}^{Y_{2}})$ be a subobject of $Y_{2}$. Let $%
x:X_{1}\rightarrow X_{2}$  and $y:Y_{1}\rightarrow Y_{2}$ be
morphisms such that $y\circ
i_{X_{1}}^{Y_{1}}=i_{X_{2}}^{Y_{2}}\circ x$. Then there exists
a unique morphism, which we denote by $y/x=\frac{y}{x}:Y_{1}/{X_{1}}%
\rightarrow Y_{2}/{X_{2},}$ such that $\frac{y}{x}\circ
p_{X_{1}}^{Y_{1}}=p_{X_{2}}^{Y_{2}}\circ y$:
\begin{equation*}
\xymatrix@R=20pt@C=40pt{
  X_1  \ar[d]_{x} \ar@{^{(}->}[r]^{i_{X_1}^{Y_1}} & Y_1 \ar[d]_{y} \ar[r]^{p_{X_1}^{Y_1}} & \frac{Y_1}{X_1} \ar[d]^{\frac{y}{x}} \\
  X_2  \ar@{^{(}->}[r]^{i_{X_2}^{Y_2}} & Y_2 \ar[r]^{p_{X_2}^{Y_2}} & \frac{Y_2}{X_2}   }
\end{equation*}

\section{Wedge products in Monoidal Categories}

\begin{claim}
\label{MonCat} Let us recall that a \emph{monoidal category} is a
category $\mathcal{M}$ that is endowed with a functor $\otimes
:\mathcal{M}\times \mathcal{M}\rightarrow \mathcal{M}$, an object
$\mathbf{1}\in \mathcal{M}$ and functorial isomorphisms:
$a_{X,Y,Z}:(X\otimes Y)\otimes Z\rightarrow X\otimes (Y\otimes
Z),$ $l_{X}:\mathbf{1}\otimes X\rightarrow X$ and $r_{X}:X\otimes \mathbf{1}%
\rightarrow X.$ The functorial morphism $a$ is called the \emph{%
associativity constraint }and\emph{\ }satisfies the \emph{Pentagon Axiom, }%
that is the following diagram
\begin{equation*}
                  \xymatrix@R=45pt@C=-30pt{
                  &((U\otimes V)\otimes W)\otimes X
                  \ar[rr]^{ \alpha_{U,V,W}\otimes X}
                  \ar[dl]|{ \alpha_{U\otimes V,W,X}}
                  &&(U\otimes (V\otimes W))\otimes X
                  \ar[dr]|{ \alpha_{U,V\otimes W,X}}&\\
                  (U\otimes V)\otimes (W\otimes X)
                  \ar[drr]|{ \alpha_{U,V,W\otimes X}}
                  &&&&U\otimes ((V\otimes W)\otimes X)
                  \ar[dll]|{ U\otimes \alpha_{V,W,X}}
                  \\&&U\otimes (V\otimes (W\otimes X))&&&
                 }
\end{equation*}
is commutative, for every $U,\,V,$ $W,$ $X$ in $\mathcal{M}.$ The morphisms $%
l$ and $r$ are called the \emph{unit constraints} and they are
assumed to satisfy the \emph{Triangle Axiom, }i.e. the following
diagram
\begin{equation*}\xymatrix@R=30pt@C=-2pt{
                 (V\ot \mathbf{1})\ot W \ar[dr]_{r_V\ot
                 W}\ar[rr]^{a_{V,\mathbf{1},W}}&&V\ot (\mathbf{1}\ot
                 W)\ar[dl]^{V\ot l_W} &\\
                 &V\ot W
                 }
\end{equation*}
is commutative. The object $\mathbf{1}$ is called the \emph{unit} of $%
\mathcal{M}$. For details on monoidal categories we refer to \cite[Chapter XI]%
{Ka} and \cite{Maj2}. A monoidal category is called \emph{strict}
if the associativity constraint and unit constraints are the
corresponding identity morphisms.
\end{claim}

\begin{claim}
\label{cl:CohThm}As it is noticed in \cite[p. 420]{Maj2}, the
Pentagon Axiom solves the consistency problem that appears because
there are two ways to go from $((U\otimes V)\otimes W)\otimes X$
to $U\otimes (V\otimes (W\otimes X)).$ The coherence theorem, due
to S. Mac Lane, solves the similar problem for the tensor product
of an arbitrary number of objects in $\mathcal{M}.$ Accordingly
with this theorem, we can always omit all brackets and simply
write $X_{1}\otimes \cdots \otimes X_{n}$ for any object obtained from $%
X_{1},\ldots ,X_{n}$ by using $\otimes $ and brackets. Also as a
consequence of the coherence theorem, the morphisms $a,$ $l,$ $r$
take care of themselves, so they can be omitted in any computation
involving morphisms in $\mathcal{M.}$\newline The notions of
algebra, module over an algebra, coalgebra and comodule over a
coalgebra can be introduced in the general setting of monoidal
categories. For more details, see \cite{AMS}.
\end{claim}

We quote from \cite[2.4]{Cotensor} the following definition. We
remark that, in this context, we prefer to use the word
"coabelian" instead of "abelian".

\begin{definition}\label{abelianmonoidal}A monoidal category
$(\M,\ot,\mathbf{1})$ will be called a \emph{\textbf{coabelian
monoidal category} } if:
\begin{enumerate}
    \item $\M$ is an abelian category
    \item both the functors $X\ot (-):\M\to\M$ and $(-)\ot
    X:\M\to\M$ are additive and left exact, for every object $X\in
    \M$.
\end{enumerate}
\end{definition}

\begin{claim}
Let $E$ be a coalgebra in a coabelian monoidal category
$\mathcal{M}$.\\ Let us recall,
(see \cite[page 60]{Mo}), the definition of wedge of two subobjects $%
X,Y $ of $E$ in $\M:$%
\begin{equation*}
(X\wedge_E Y,i_{X\w Y}^E):=Ker[ (p _{X}\otimes p _{Y}) \circ
\triangle _{E}] ,
\end{equation*}
where $p _{X}:E\rightarrow E/X$ and $p _{Y}:E\rightarrow E/Y$ are
the canonical quotient maps. In particular we have the following
exact sequence:
\begin{equation*}
\xymatrix@C=0.9cm{
  0 \ar[r] & X\wedge_E Y \ar[rr]^{i_{X\w Y}^E} && E \ar[rr]^(.4){(p _{X}\otimes p _{Y}) \circ \triangle
  _{E}} && E/X \ot E/Y.}
\end{equation*}
Consider the following commutative diagrams in $\M$
  \begin{equation*}
\begin{tabular}{cc}
$\xymatrix@C=2cm{
  X_1 \ar[d]_{x} \ar@{^{(}->}[r]^{i_{X_1}^{E_1}}
                & E_1 \ar[d]^{e}  \\
  X_2 \ar@{^{(}->}[r]_{i_{X_2}^{E_2}}
                & E_2             }
$& $ \xymatrix@C=2cm{
  Y_1 \ar[d]_{y} \ar@{^{(}->}[r]^{i_{Y_1}^{E_1}}
                & E_1 \ar[d]^{e}  \\
  Y_2 \ar@{^{(}->}[r]_{i_{Y_2}^{E_2}}
                & E_2             }
$
\end{tabular}
\end{equation*} where $e$ is a coalgebra homomorphism. Then there is a unique
morphism $x\wedge_e y:X_1\wedge_{E_1} Y_1\to X_2\wedge_{E_2} Y_2$
such that the following diagram
\begin{equation*}
\xymatrix@C=2cm{
  X_1\wedge_{E_1} Y_1 \ar@{.>}[d]_{x\wedge_e y} \ar[r]^{i_{X_1\wedge_{E_1} Y_1}^{E_1}}
                & E_1 \ar[d]^{e}  \\
  X_2\wedge_{E_2} Y_2 \ar[r]_{i_{X_2\wedge_{E_2} Y_2}^{E_2}}
                & E_2             }
\end{equation*}
commutes. In fact we have
\begin{eqnarray*}
&&(p_{X_{2}}^{E_{2}}\otimes p_{Y_{2}}^{E_{2}})\circ \Delta
_{E_{2}}\circ
e\circ i_{X_1\wedge_{E_1} Y_1}^{E_1} \\
&=&(p_{X_{2}}^{E_{2}}\otimes p_{Y_{2}}^{E_{2}})\circ (e\otimes
e)\circ \Delta
_{E_{1}}\circ i_{X_1\wedge_{E_1} Y_1}^{E_1} \\
&=&(\frac{e}{x}\otimes \frac{e}{y})\circ (p_{X_{1}}^{E_{1}}\otimes
p_{Y_{1}}^{E_{1}})\circ \Delta _{E_{1}}\circ i_{X_1\wedge_{E_1}
Y_1}^{E_1}=0
\end{eqnarray*}
so that, since $(X_2\wedge_{E_2} Y_2,i_{X_2\wedge_{E_2}
Y_2}^{E_2})$ is the kernel of $(p_{X_{2}}^{E_{2}}\otimes
p_{Y_{2}}^{E_{2}})\circ \Delta _{E_{2}}$, we conclude.
\end{claim}

\begin{lemma}
Consider the following commutative diagrams in $\M$
  \begin{equation*}
\begin{tabular}{cc}
$\xymatrix@C=2cm{
  X_1 \ar[d]_{x} \ar@{^{(}->}[r]^{i_{X_1}^{E_1}}
                & E_1 \ar[d]^{e}  \\
  X_2 \ar[d]_{x'}\ar@{^{(}->}[r]_{i_{X_2}^{E_2}}
                & E_2 \ar[d]_{e'}           \\
  X_3 \ar@{^{(}->}[r]_{i_{X_3}^{E_3}}
                & E_3 }
$& $\xymatrix@C=2cm{
  Y_1 \ar[d]_{y} \ar@{^{(}->}[r]^{i_{Y_1}^{E_1}}
                & E_1 \ar[d]^{e}  \\
  Y_2 \ar[d]_{y'}\ar@{^{(}->}[r]_{i_{Y_2}^{E_2}}
                & E_2 \ar[d]_{e'}           \\
  Y_3 \ar@{^{(}->}[r]_{i_{Y_3}^{E_3}}
                & E_3 }$
\end{tabular}
\end{equation*}  where $e$ and $e'$ are coalgebra
homomorphisms. Then we have
\begin{equation}\label{formula compos of wedge}
  (x'\wedge_{e'} y')\circ  (x\wedge_{e} y)= (x'x\wedge_{e'e} y'y)
\end{equation}
\end{lemma}

\begin{proof}: straightforward.
\end{proof}

We now recall some definitions and some (standard) results
established in \cite{Cotensor}.

\begin{claim}
\label{def of delta_n}Let $X$ be an object in a coabelian monoidal category $%
\left( \mathcal{M},\otimes ,\mathbf{1}\right) $. Set%
\begin{equation*}
X^{\otimes 0}=\mathbf{1},\qquad X^{\otimes 1}=X\qquad
\text{and}\qquad X^{\otimes n}=X^{\otimes n-1}\otimes X,\text{ for
every }n>1
\end{equation*}%
and for every morphism $f:X\rightarrow Y$ in $\mathcal{M}$, set%
\begin{equation*}
f^{\otimes 0}=\mathrm{Id}_{\mathbf{1}},\qquad f^{\otimes 1}=f\qquad \text{and%
}\qquad f^{\otimes n}=f^{\otimes n-1}\otimes f,\text{ for every
}n>1.
\end{equation*}%
Let $\left( C,\Delta _{C},\varepsilon _{C}\right) $ be a coalgebra in $%
\mathcal{M}$ and for every $n\in
\mathbb{N}
,$ define the $n^{\text{th}}$ iterated comultiplication of $C,$
$$\Delta _{C}^{n}:C\rightarrow C^{\otimes {n+1}},$$ by
\begin{equation*}
\Delta _{C}^{0}=\text{Id}_{C},\qquad \Delta _{C}^{1}=\Delta _{C}\qquad \text{%
and}\qquad \Delta _{C}^{n}=\left( \Delta _{C}^{\otimes n-1}\otimes
C\right) \Delta _{C},\text{ for every }n>1.
\end{equation*}%
Let $\delta :D\rightarrow E$ be a monomorphism which is a
homomorphism of
coalgebras in $\mathcal{M}$. Denote by $(L,p)$ the cokernel of $\delta $ in $%
\mathcal{M}$. Regard $D$ as a $E$-bicomodule via $\delta $ and
observe that $L$ is a $E$-bicomodule and $p$ is a morphism of
bicomodules. Let
\begin{equation*}
(D^{\wedge _{E}^{n}},\delta _{n}):=\ker (p^{\otimes {n}}\Delta
_{E}^{n-1})
\end{equation*}%
for any $n\in \mathbb{N}\setminus \{0\}.$ Note that $(D^{\wedge
_{E}^{1}},\delta _{1})=(D,\delta )$ and $(D^{\wedge
_{E}^{2}},\delta _{2})=D\wedge _{E}D.$ \newline In order to
simplify the notations we set $(D^{\wedge _{E}^{0}},\delta
_{0})=(0,0).$\newline Now, since $\mathcal{M}$ has left exact
tensor functors and since $p^{\otimes {n}}\Delta _{E}^{n-1}$ is a
morphism of $E$-bicomodules (as a composition of morphisms
of $E$-bicomodules), we get that $D^{\wedge _{E}^{n}}$ is a coalgebra and $%
\delta _{n}:D^{\wedge _{E}^{n}}\rightarrow E$ is a coalgebra
homomorphism for any $n>0$ and hence for any $n\in \mathbb{N}$.
\end{claim}

\begin{proposition}\cite[Proposition 1.10]{Cotensor}\label{pro: limit of delta} Let $\delta:D\to E$ be a monomorphism which is a morphism of coalgebras in a coabelian monoidal category
$\M$. Then, for any $i\leq j$ in $\mathbb{N}$, there is a (unique)
morphism $\xi_{i}^j:D^{\wedge_E ^i}\to D^{\wedge_E ^j}$ such that
\begin{equation}\label{compatibility of delta i}
\delta_j\xi_{i}^j=\delta_i.
\end{equation} Moreover $\xi_{i}^j$ is a coalgebra
homomorphism and $((D^{\wedge_E ^i})_{i\in \mathbb{N}},(\xi
_{i}^j)_{i,j\in \mathbb{N}})$ is a direct system in $\M$ whose
direct limit, if it exists, carries a natural coalgebra structure
that makes it the direct limit of $((D^{\wedge_E ^i})_{i\in
\mathbb{N}},(\xi _{i}^j)_{i,j\in \mathbb{N}})$ as a direct system
of coalgebras.
\end{proposition}

\begin{proposition}\cite[Proposition
2.17]{Cotensor} \label{pro: D^2}Let $\delta :D\rightarrow E$ be a
monomorphism which is a coalgebra homomorphism in a coabelian
monoidal category $\mathcal{M}$. Then we have
\begin{equation}
(D^{\wedge _{E}^{m}}\wedge _{E}D^{\wedge _{E}^{n}},i_{D^{\wedge
_{E}^{m}}\wedge _{E}D^{\wedge _{E}^{n}}}^E)=(D^{\wedge
_{E}^{m+n}},i^E_{D^{\wedge _{E}^{m+n}}}). \label{formula 2 pro:
D^2}
\end{equation}
\end{proposition}

\begin{notation}\label{notation tilde}Let $\delta:D\to E$ be a morphism of coalgebras in a cocomplete coabelian monoidal category
$\M$ with left exact tensor functors . By Proposition \ref{pro:
limit of delta} $((D^{\wedge_E ^i})_{i\in \mathbb{N}},(\xi
_{i}^j)_{i,j\in \mathbb{N}})$ is a direct system in $\M$ whose
direct limit carries a natural coalgebra structure that makes it
the direct limit of $((D^{\wedge_E ^i})_{i\in \mathbb{N}},(\xi
_{i}^j)_{i,j\in \mathbb{N}})$ as a direct system of
coalgebras.\\From now on we will use the following notation
\begin{eqnarray*}
D^{n}:=D^{\w n},\text{ for every }n\in \mathbb{N},\\
(\widetilde{D}_E, (\xi_i)_{i\in \mathbb{N}}):=
\underrightarrow{\lim }(D^{\wedge_E ^i})_{i\in
\mathbb{N}},\end{eqnarray*}
 where
$\xi_i:D^{\wedge_E ^i}\to \widetilde{D}_E$ denotes the structural
morphism of the direct limit. We note that, since
$\widetilde{D}_E$ is a direct limit of coalgebras, the canonical
(coalgebra) homomorphisms $(\delta _{i}:D^{\wedge_E ^i}\rightarrow
E)_{i\in \mathbb{N}}$, which are compatible by (\ref{compatibility
of delta i}), factorize to a unique coalgebra homomorphism
$$\widetilde{\delta} :\widetilde{D}_E\rightarrow E$$ such that
$\widetilde{\delta} \xi_i=\delta_i$, for any $i\in \mathbb{N}$.
\end{notation}

\section{The Heyneman-Radford Theorem for Monoidal
Categories}\label{section: Montgomery}

\begin{definition}
Let $E$ be a coalgebra and let $\delta :X\rightarrow E$ be a
monomorphism in a coabelian monoidal category $\mathcal{M}$.
Define the morphism
\begin{equation*}
\alpha _{X}^{E}:E\rightarrow \frac{E}{X}\otimes \frac{E}{X}
\end{equation*}%
by setting
\begin{equation*}
\alpha _{X}^{E}=\left( p_{X}^{E}\otimes p_{X}^{E}\right) \circ
\Delta _{E}.
\end{equation*} Observe that $(X\wedge_E X,i_{X\wedge _E X}^E) =\K(\alpha_X^E).$
\end{definition}

\begin{lemma}
\label{lem: alpha}Let $\delta :D\rightarrow E$ and let
$f:E\rightarrow C$ be coalgebra homomorphisms in a coabelian
monoidal category $\mathcal{M}$. Assume that both $\delta $ and
$f\circ \delta $ are monomorphism. Then the following diagram
\begin{equation*}
\xymatrix@R=35pt@C=50pt{
  E \ar[d]_{\alpha_D^E} \ar[r]^{f} & C \ar[d]^{\alpha_D^C} \\
  \frac{E}{D}\otimes\frac{E}{D} \ar[r]^{\frac{f}{D}\otimes\frac{f}{D}} & \frac{C}{D}\otimes\frac{C}{D}   }
\end{equation*}%
is commutative.
\end{lemma}

\begin{proof}
Note that the notations $E/D$ and $C/D$ make sense as both $\delta $ and $%
f\circ \delta $ are monomorphisms. We have%
\begin{eqnarray*}
\left( \frac{f}{D}\otimes \frac{f}{D}\right) \circ \alpha _{D}^{E}
&=&\left( \frac{f}{D}\otimes \frac{f}{D}\right) \circ \left(
p_{D}^{E}\otimes
p_{D}^{E}\right) \circ \Delta _{E} \\
&=&\left( p_{D}^{C}\otimes p_{D}^{C}\right) \circ \left( f\otimes
f\right)
\circ \Delta _{E} \\
&=&\left( p_{D}^{C}\otimes p_{D}^{C}\right) \circ \Delta _{C}\circ
f=\alpha _{D}^{C}\circ f.
\end{eqnarray*}
\end{proof}

\begin{lemma}
\label{lem: tau}Let $D$ and $E$ be coalgebras in a coabelian
monoidal category $\mathcal{M}$. Let $\delta :D\rightarrow E$ be a
monomorphism which is a morphism of coalgebras in $\mathcal{M}$.
Then, for every $n\in
\mathbb{N}
$, there exists a unique morphism $\tau _{n}:D^{n+1}\rightarrow
D^{n}/D\otimes D^{n}/D$ such that the following diagram
\begin{equation*}
\xymatrix@R=40pt@C=80pt{
                &         D^{n+1}\ar@{.>}[dl]_{\tau_n} \ar[d]^{\alpha_D^{D^{n+1}}}     \\
  \frac{D^{n}}{D}\otimes\frac{D^{n}}{D} \ar[r]_{\frac{\xi_n^{n+1}}{D}\otimes%
\frac{\xi_n^{n+1}}{D}} & \frac{D^{n+1}}{D}\otimes\frac{D^{n+1}}{D}
} \end{equation*}%
is commutative.
\end{lemma}

\begin{proof}
Consider the following exact sequence%
\begin{equation}\label{formula:
exact1} \xymatrix@C=0.5cm{
  0 \ar[r] & \frac{D^{n}}{D} \ar[rr]^{\frac{\xi
_{n}^{n+1}}{D}} && \frac{D^{n+1}}{D} \ar[rr]^{\frac{D^{n+1}}{\xi
_{1}^{n}}} && \frac{D^{n+1}}{D^{n}} \ar[r] & 0 }
\end{equation}
By applying the functor $D^{n+1}/D\otimes \left( -\right) $ we get%
\begin{equation*}
\xymatrix@C=0.9cm{
  0 \ar[r] & \frac{D^{n+1}}{D}\otimes \frac{D^{n}}{D} \ar[rr]^{\frac{D^{n+1}}{%
D}\otimes \frac{\xi _{n}^{n+1}}{D}} && \frac{D^{n+1}}{D}%
\otimes \frac{D^{n+1}}{D} \ar[rr]^{\frac{D^{n+1}}{D}\otimes \frac{D^{n+1}}{%
\xi _{1}^{n}}} && \frac{D^{n+1}}{D}\otimes \frac{D^{n+1}}{D^{n}}
\ar[r] & 0 }
\end{equation*}
We have%
\begin{eqnarray*}
&&\left( \frac{\delta _{n+1}}{D}\otimes \frac{\delta
_{n+1}}{D^{n}}\right) \circ \left( \frac{D^{n+1}}{D}\otimes
\frac{D^{n+1}}{\xi _{1}^{n}}\right)
\circ \alpha _{D}^{D^{n+1}} \\
&=&\left( \frac{\delta _{n+1}}{D}\otimes \frac{\delta
_{n+1}}{D^{n}}\right) \circ \left( \frac{D^{n+1}}{D}\otimes
\frac{D^{n+1}}{\xi _{1}^{n}}\right) \circ \left(
p_{D}^{D^{n+1}}\otimes p_{D}^{D^{n+1}}\right) \circ \Delta
_{D^{n+1}} \\
&=&\left( \frac{\delta _{n+1}}{D}\otimes \frac{\delta
_{n+1}}{D^{n}}\right) \circ \left( p_{D}^{D^{n+1}}\otimes
p_{D^{n}}^{D^{n+1}}\right) \circ \Delta
_{D^{n+1}} \\
&=&\left( p_{D}^{E}\otimes p_{D^{n}}^{E}\right) \circ \left(
\delta
_{n+1}\otimes \delta _{n+1}\right) \circ \Delta _{D^{n+1}} \\
&=&\left( p_{D}^{E}\otimes p_{D^{n}}^{E}\right) \circ \Delta
_{E}\circ \delta _{n+1}=0.
\end{eqnarray*}%
In fact, by Proposition \ref{pro: D^2}, $D^{n+1}=D\wedge
_{E}D^{n}.$ Since $\frac{\delta _{n+1}}{D}\otimes
\frac{\delta _{n+1}}{D^{n}}$ is a monomorphism, we obtain%
\begin{equation}
\left( \frac{D^{n+1}}{D}\otimes \frac{D^{n+1}}{\xi
_{1}^{n}}\right) \circ \alpha _{D}^{D^{n+1}}=0  \label{formula: xi
and alpha}
\end{equation}%
so that, as the above sequence is exact, by the universal property
of kernels, there exists a unique morphism
\begin{equation*}
\beta _{n}:D^{n+1}\rightarrow \frac{D^{n+1}}{D}\otimes
\frac{D^{n}}{D}
\end{equation*}%
such that%
\begin{equation}
\left( \frac{D^{n+1}}{D}\otimes \frac{\xi _{n}^{n+1}}{D}\right)
\circ \beta _{n}=\alpha _{D}^{D^{n+1}}.  \label{formula: def beta
n}
\end{equation}%
By applying the functor $\left( -\right) \otimes D^{n}/D$ to
(\ref{formula: exact1}), we get
\begin{equation*}
\xymatrix@C=0.9cm{
  0 \ar[r] & \frac{D^{n}}{D}\otimes
\frac{D^{n}}{D} \ar[rr]^{\frac{\xi _{n}^{n+1}}{D}\otimes
\frac{D^{n}}{D}} && \frac{D^{n+1}}{D}%
\otimes \frac{D^{n}}{D} \ar[rr]^{\frac{D^{n+1}}{\xi _{1}^{n}}\otimes \frac{%
D^{n}}{D}} &&\frac{D^{n+1}}{D^{n}}\otimes \frac{D^{n}}{D} \ar[r] &
0 .}\end{equation*}
We have%
\begin{eqnarray*}
&&\left( \frac{D^{n+1}}{D^{n}}\otimes \frac{\xi
_{n}^{n+1}}{D}\right) \circ \left( \frac{D^{n+1}}{\xi
_{1}^{n}}\otimes \frac{D^{n}}{D}\right) \circ
\beta _{n} \\
&=&\left( \frac{D^{n+1}}{\xi _{1}^{n}}\otimes
\frac{D^{n+1}}{D}\right) \circ \left( \frac{D^{n+1}}{D}\otimes
\frac{\xi _{n}^{n+1}}{D}\right) \circ \beta
_{n} \\
&\overset{(\ref{formula: def beta n})}{=}&\left( \frac{D^{n+1}}{\xi _{1}^{n}}%
\otimes \frac{D^{n+1}}{D}\right) \circ \alpha _{D}^{D^{n+1}} =0
\end{eqnarray*}%
where the last equality can be proved similarly to (\ref{formula:
xi and alpha}). Since $\frac{D^{n+1}}{D^{n}}\otimes \frac{\xi
_{n}^{n+1}}{D}$ is a
monomorphism we get%
\begin{equation*}
\left( \frac{D^{n+1}}{\xi _{1}^{n}}\otimes \frac{D^{n}}{D}\right)
\circ \beta _{n}=0
\end{equation*}%
so that, as the previous sequence is exact, by the universal
property of kernels there exists a unique morphism
\begin{equation*}
\tau _{n}:D^{n+1}\rightarrow \frac{D^{n}}{D}\otimes
\frac{D^{n}}{D}
\end{equation*}%
such that%
\begin{equation*}
\left( \frac{\xi _{n}^{n+1}}{D}\otimes \frac{D^{n}}{D}\right)
\circ \tau _{n}=\beta _{n}.
\end{equation*}%
Finally we have%
\begin{eqnarray*}
\left( \frac{\xi _{n}^{n+1}}{D}\otimes \frac{\xi
_{n}^{n+1}}{D}\right) \circ \tau _{n} &=&\left(
\frac{D^{n+1}}{D}\otimes \frac{\xi _{n}^{n+1}}{D}\right) \circ
\left( \frac{\xi _{n}^{n+1}}{D}\otimes \frac{D^{n}}{D}\right)
\circ
\tau _{n} \\
&=&\left( \frac{D^{n+1}}{D}\otimes \frac{\xi _{n}^{n+1}}{D}\right)
\circ \beta _{n}=\alpha _{D}^{D^{n+1}}.
\end{eqnarray*}
\end{proof}

\begin{theorem}
\label{teo: Montgomery}Let $D$ and $E$ be coalgebras in a
cocomplete coabelian monoidal category $\mathcal{M}$ satisfying
AB5. Let $\delta :D\rightarrow E$ be a monomorphism which is a
morphism of coalgebras in $\mathcal{M}$ and keep the notations
introduced in Notation \ref{notation tilde}.\\ Let $f:E\rightarrow
C$ be a coalgebra homomorphism and assume that
\begin{equation*}
f\circ \delta _{2}:D\wedge _{E}D\rightarrow C
\end{equation*}%
is a monomorphism. Then the coalgebra homomorphism%
\begin{equation*}
f\circ \widetilde{\delta }:\widetilde{D}_{E}\rightarrow C
\end{equation*}%
is a monomorphism.
\end{theorem}

\begin{proof}
Since $\mathcal{M}$ satisfies AB5, it is enough to prove that
$f\circ \widetilde{\delta }\circ \xi _{n}=f\circ \delta _{n}$ is a
monomorphism for every $n\in
\mathbb{N}
.$

For $n=0$, we have $f\circ \delta _{0}=f\circ 0=0$ which is a
monomorphism as $D^{0}=0$.

For $n=1$, we have $f\circ \delta _{1}=f\circ \delta _{2}\circ \xi
_{1}^{2}$ which is a monomorphism.

Let $n\geq 2$ and let us assume that $f\circ \delta _{n}$ is a
monomorphism. Let us prove that $f\circ \delta _{n+1}$ is a
monomorphism. Let $\lambda
:X\rightarrow D^{n+1}$ be a morphism such that%
\begin{equation*}
f\circ \delta _{n+1}\circ \lambda =0
\end{equation*}%
and consider the following diagram%
\begin{equation*}
\xymatrix@R=50pt@C=60pt{
             & D\wedge_{D^{n+1}}D \ar[d]^{i_{D\wedge_{D^{n+1}}D}^{D^{n+1}}} \ar[r]^{D\wedge_{\delta_{n+1}}D} & D\wedge_{E}D \ar[d]^{\delta_2} \\
  X \ar@{.>}[ru]^{\overline{\lambda}} \ar[r]^{\lambda} & D^{n+1} \ar[dl]_{\tau_n}\ar[d]^{\alpha_D^{D^{n+1}}}
  \ar[r]^{\delta_{n+1}} & E \ar[r]^{f}\ar[d]^{\alpha_D^E} &C\ar[d]^{\alpha_D^C}\\
  \frac{D^{n}}{D}\otimes\frac{D^{n}}{D} \ar[r]_{\frac{\xi_n^{n+1}}{D}\otimes%
\frac{\xi_n^{n+1}}{D}} & \frac{D^{n+1}}{D}\otimes\frac{D^{n+1}}{D}
\ar[r]_{\frac{\delta_{n+1}}{D}\otimes\frac{\delta_{n+1}}{D}} & \frac{E}{D}\otimes%
\frac{E}{D} \ar[r]_{\frac{f}{D}\otimes\frac{f}{D}} & \frac{C}{D}\otimes%
\frac{C}{D}   }
\end{equation*}
where all the squares are commutative in view of Lemma \ref{lem:
alpha} and
the bottom triangle commutes in view of Lemma \ref{lem: tau}. We have%
\begin{eqnarray*}
&&\left( \frac{f\delta _{n}}{D}\otimes \frac{f\delta
_{n}}{D}\right) \circ
\tau_n \circ \lambda  \\
&=&\left( \frac{f\delta _{n+1}\xi _{n}^{n+1}}{D}\otimes
\frac{f\delta
_{n+1}\xi _{n}^{n+1}}{D}\right) \circ \tau_n \circ \lambda  \\
&=&\left( \frac{f}{D}\otimes \frac{f}{D}\right) \circ \left(
\frac{\delta _{n+1}}{D}\otimes \frac{\delta _{n+1}}{D}\right)
\circ \left( \frac{\xi _{n}^{n+1}}{D}\otimes \frac{\xi
_{n}^{n+1}}{D}\right) \circ \tau_n \circ
\lambda  \\
&=&\alpha _{D}^{C}\circ f\circ \delta _{n+1}\circ \lambda =0.
\end{eqnarray*}%
Since $f\circ \delta _{n}$ is a monomorphism, we get that also
$f\delta
_{n}/D\otimes f\delta _{n}/D$ is a monomorphism so that we obtain%
\begin{equation*}
\tau_n \circ \lambda =0
\end{equation*}%
and we have%
\begin{equation*}
\alpha _{D}^{D^{n+1}}\circ \lambda =\left( \frac{\xi
_{n}^{n+1}}{D}\otimes \frac{\xi _{n}^{n+1}}{D}\right) \circ \tau_n
\circ \lambda =0.
\end{equation*}%
Thus, since $\left( D\wedge _{D^{n+1}}D,i_{D\wedge
_{D^{n+1}}D}^{D^{n+1}}\right) =\ker \left( \alpha
_{D}^{D^{n+1}}\right) $,
by the universal property of the kernel, there exists a unique morphis, $%
\overline{\lambda }:X\rightarrow D\wedge _{D^{n+1}}D$ such that%
\begin{equation*}
\lambda =i_{D\wedge _{D^{n+1}}D}^{D^{n+1}}\circ \overline{\lambda
}.
\end{equation*}%
Now we have%
\begin{equation*}
f\circ \delta _{2}\circ \left( D\wedge _{\delta _{n+1}}D\right)
\circ \overline{\lambda }=f\circ \delta _{n+1}\circ \lambda =0.
\end{equation*}%
Since $f\circ \delta _{2}$ and $D\wedge _{\delta _{n+1}}D\ $are
monomorphisms, we get that $\overline{\lambda }=0$ and hence
$\lambda =0$.
\end{proof}

\begin{corollary}(Heyneman-Radford) (\cite[Proposition 2.4.2]{HR} or \cite[Theorem 5.3.1, page
65]{Mo})
Let $K$ be a field. Let $E$ and $C$ be $K$-coalgebras and let $%
f:E\rightarrow C$ be a coalgebra homomorphism such that $f_{\mid
D\wedge _{E}D}$ is injective, where $D$ is the coradical of $E$.
Then $f$ is injective.
\end{corollary}

\begin{proof}
Since $D$ is the coradical of $E$ is well known that
$(E,\Id_E)=(\widetilde{D}_E,\widetilde{\delta)}$ (see e.g.
\cite[Corollary 9.0.4, page 185]{Sw}). The conclusion follows by
Theorem \ref{teo: Montgomery} applied in the case when
$\mathcal{M}$ is the category of vector spaces over $K.$ Observe
that in this case "monomorphism" is equivalent to "injective".
\end{proof}

\end{document}